\documentclass[a4paper,11pt]{article}
\usepackage{charter}
\usepackage{latexsym}
\usepackage{amsmath}
\usepackage{amsfonts}
\usepackage{amssymb}
\usepackage{vmargin}
\usepackage{url}
\setmarginsrb{2cm}{2cm}{2cm}{2cm}{0mm}{0mm}{0mm}{10mm}
\usepackage{authblk}
\usepackage[round]{natbib}
\usepackage{color}
\usepackage[colorlinks=true, allcolors=blue]{hyperref}

\begin{document}
\title{\Large \bf Challenges in teaching Real Analysis classes at\\ the University of PGRI, South Sumatra, Indonesia}
\author[1]{E. Septiati}
\author[2]{N. Karjanto\thanks{\url{natanael.karjanto@nottingham.edu.my}}}
\affil[1]{\small Department of Mathematics Education, University of PGRI, Palembang, South Sumatra, Indonesia}
\affil[2]{\small School of Applied Mathematics, Faculty of Engineering, The University of Nottingham Malaysia Campus Semenyih, Selangor, Malaysia}
\date{\footnotesize Presented July 25, 2008; \ Submitted July 30, 2008; \  Updated \today}
\maketitle

\begin{abstract}
This paper discusses our experiences and challenges in teaching advanced undergraduate Real Analysis classes for Mathematics Education students at the University of PGRI (\textsl{Persatuan Guru Republik Indonesia}, Indonesian Teachers Association) Palembang, South Sumatra, Indonesia. We observe that the syllabus contains topics with a high level of difficulty for the students who are specialized in education and intend to teach mathematics at the secondary level. The conventional lecturing method is mainly implemented during the class, with some possible variations of the method, including the Texas method (also known as Moore's method) and the small group guided discovery method. In particular, the latter method has been implemented successfully for a Real Analysis class at Dartmouth College, New Hampshire by~\citeauthor{Dumitrascu09} in 2006. Although it is a real challenge to apply a specific teaching method that will be able to accommodate a large number of students, the existing teaching activities can still be improved and a more effective method could be implemented in the future. Furthermore, the curriculum contents should be adapted for an audience in Mathematics Education to equip them for their future career as mathematics teachers. Any constructive suggestions are welcome for the improvement of our mathematics education system at the university as well as on the national scale.
\end{abstract}

\section{Introduction}

There are several teaching methods and certainly teaching style varies from one lecturer to another. Several ways of teaching methods that are commonly carried out in many parts of the world include the following: questioning, explaining, demonstrating, collaborating, and learning by teaching~\citep{Committee97,Good08}. In particular, the learning by teaching (German, \emph{Lernen durch Lehren}) is a widespread method in Germany, where the students take the teacher's role and teach their peers.

More specifically, we want to use certain methods in teaching mathematics, and in this context, in teaching Real Analysis classes. Methods of teaching mathematics include the following: classical education, rote learning, exercises, problem-solving, new math, historical method, and reform or standard-based mathematics~\citep{Clarke03, Fan04, Lockhart09}. In particular, there are  significant research results on the implementation of the realistic mathematics education method in Indonesia~\citep{Armanto02, Fauzan02, Hadi02, Zulkardi02}. Furthermore, cooperative learning methods are now being used more and more often in teaching undergraduate mathematics and science~\citep{Davidson91, Rogers01, Dubinsky97, Finkel83} as well as in higher education settings~\citep{Ledlow99, Milis10}.

Another teaching method which accompanies a conventional instruction method is known as the `guided discovery' method. In this method, the students learn through personal experience\,\dots\,with limited guidance from the lecturer,\,\dots\,thought-provoking topics are introduced as questions for investigation by the students~\citep{Davidson70, Davidson71}. Experience of implementing teaching methods using a combination of guided discovery, lecturing, and group work in an undergraduate Real Analysis class has been proven to improve students' understanding~\citep{Dumitrascu99}. The author concludes that the guided discovery method is an excellent modality of exposing students to mathematical research.

We observe that a conventional teaching method using instruction and lecturing for the Real Analysis courses presents a challenge for the students who are specialized in Mathematics Education. This challenge motivates us in bringing this problem into the surface.

What is ``real analysis'' and what is the scope of the course on Real Analysis? Real analysis is a branch of mathematical analysis dealing with the set of real numbers. In particular, it deals with the analytic properties of real functions and sequences, including convergence and limits of sequences of real numbers, the calculus of the real numbers and continuity, smoothness, and related properties of real-valued functions~\citep{Bressoud07, Krantz04, Stein09, Trench13}. Certainly, a course on Real Analysis should cover the aforementioned materials. This course is an important component of mathematics curriculum for both educational and noneducational streams at the undergraduate level.

In this paper, we share our experiences in implementing different teaching methods to the Real Analysis courses. This paper is organized as follows. Section~\ref{org} discusses the organization of the courses, including reference textbooks being used and the method of assessment. Section~\ref{pro} explains the challenges and difficulties that students face in following the classes. Furthermore, Section~\ref{obs} discusses our observations in conducting the classes and implementing several teaching methods. This section also provides the students' responses toward different teaching methods. Finally, Section~\ref{con} gives the conclusion and remark for future research to our discussion.

\section{Course organization} \label{org}

The classes of Real Analysis I and II are compulsory subjects for undergraduate students in Mathematics Education at the University of PGRI, Palembang, South Sumatra, Indonesia. These courses carry three credit points and are given to third-year students or in the fifth and sixth semesters of their study. There is only one-time interaction every week and it lasts for 150 minutes, which is three times 50 minutes.

There are ten topics in total which are covered in the Real Analysis courses, in which six topics belong to Real Analysis I and the other four belong to Real Analysis II. The materials covered in Real Analysis I are ordered set, field, Euclidean space, metric space, topological concepts in metric space, and sets in metric space. The topics covered in Real Analysis II are convergence sequence, Darboux integral, Riemann integral, and Rieman-Stieltjes integral.

The number of students in one class ranges from 15 up to 40 students and there are 9 parallel classes for the same course with a total number of 335 students. One senior lecturer plays a role as a coordinator for three other lecturers. Practically, the material presented in this paper is merely based on our observation of four different classes taught by one of us (ES). One class consists of only 15 students while the other three are 28, 29, and 32 students, respectively.

Lecture notes are prepared and compiled from several mathematical analysis books, among others are \textit{Introduction to Real Analysis} by~\cite{Bartle11}, \textit{Principles of Mathematical Analysis} by~\cite{Rudin76}, \textit{Pengantar Analisis Real} by~\cite{Darmwawijaya86}, and lecture notes on Real Analysis from Malang State University, East Java, Indonesia. For Real Analysis II, additional references have been used, among others are \textit{Analisis Real} by~\cite{Soemantri93}, \textit{Real Analysis} by~\cite{Royden88} and \textit{Fundamental Concepts of Analysis} by~\cite{Smith81}. There are other excellent textbooks on Real Analysis at the introductory level, including~\citep{Browder12, Kolmogorov75, Protter12, Stromberg15, Wheeden15}.

The method of assessment is based on several components with different weights. One coursework is assigned and is graded on an individual basis, this assignment carries 20 percent of the final grade. Two examinations--the mid-semester and final exams--carry 30 and 50 percents of the final grade, respectively.

\section{Students' difficulty} \label{pro}

Many students have a wrong interpretation of what mathematics subjects involve. They generally associate mathematics with counting, calculation, and computation which in turn restrict the discipline into only arithmetic. One online encyclopedia defines mathematics as the body of knowledge centered on such concepts as quantity, structure, space, and change and also the academic discipline that studies them~\cite{wiki_math}. That is why mathematics includes the use of abstraction and logical reasoning which involves rigorous deduction from appropriately chosen axioms and definitions.

We observe that the students are less familiar with theorems and how to prove them. Implementing mathematics symbols and terminologies is far from familiar. The concept of set theory is still weakly comprehended. As an example, many students are not able to distinguish simple notations such as $(a,b)$, $[a,b)$, $(a,b]$ and $\{a,b\}$. In particular, the students consider proving theorem, convergence, and Riemann integral as the most difficult topics. Bear in mind that the students who specialize in Mathematics Education spend only merely of 55\% of the total credit points on Mathematics courses for the entire study period, i.e. 84 out of 154 credit points are Mathematics courses.

Furthermore, the students also face difficulties in some technical issues, in particular, to find the literature. It is rather difficult to obtain reference books since the library has a limited amount of these books while the number of students is quite massive. The price of these textbooks is considered very expensive for all of the students. Even if the students possess textbooks, yet, since English is not the mother tongue of students, the language barrier may present another significant challenge in understanding the material.

\section{Observation and students' response} \label{obs}

We have implemented several teaching methods in conducting the Real Analysis classes. These are the conventional instruction method, Moore's method, and the guided discovery method. The conventional instruction method is implemented to the majority of the class sessions, in particular, to explain definitions and new concepts. Moore's method is implemented when explaining the properties of integral. For instance, this method is used to prove the following theorem. If $f$ is a bounded function and Darboux integrable on an interval $[a,b]$, show that $f^{2}$ is also Darboux integrable on the same interval. The guided discovery method is implemented in some theorem-proving sessions.

We observe that the students prefer the conventional instruction method more than the other two methods. The students are not able to follow Moore's method at all since none of them can answer or to give an idea in solving the theorem above. It is observed that a small number of the students could follow the guided discovery method, i.e. less than 20\%. Regarding the preference of teaching method, 34.6\% of the students prefer the conventional instruction method where the teacher only lecturing, 32.7\% prefer a variation in teaching method, 13.5\% prefer the conventional instruction and discussion, 5.8\% prefer the guided discovery method, 5.8\% prefer the conventional instruction and problem-solving and 1.9\% each for preference in question and answer session, task assignment, discussion, and self-study.

Regarding the implementation of different teaching methods in understanding the material, 48\% of the students considers it helpful, 21.2\% also considers it helpful but prefers only the conventional instruction method, 5.8\% says it can be helpful but without the discussion session, 9.6\% considers it is not helpful at all and 15.4\% says it depends on the material being covered. Furthermore, we would like to know what kind of comprehension the students acquire after completing the courses on Real Analysis. Almost 40\% of the students (39.4\%) acquires logical reasoning and improvement in theorem proving, 30.3\% acquires knowledge on the topic of integral, 18.2\% improves their understanding in mathematical symbols and 12.12\% claims do not improve at all.

Regarding the material delivery by the lecturer, almost 60\% says that it is very easy to understand (59.6\%), 32.7\% respond that it is sufficiently easy to understand, and only 7.7\% say that it is difficult to understand. Regarding the availability of the textbooks, an excellent number of 94.2\% claim that it is very helpful, 3.8\% say that it is helpful but they need some other additional references and 9.6\% say that it is not helpful. The following section gives a conclusion to our discussion.

\section{Conclusion} \label{con}

We have discussed that the Real Analysis courses are very important components in the curriculum of the Mathematics Education program. Nevertheless, a majority of the students consider that these courses are very tough and challenging. We have implemented different teaching methods to help the students to get a better understanding of the materials. Even though Moore's method and the guided discovery method have been implemented successfully in some mathematics courses in several colleges in the US, we observe that these methods are still difficult to be implemented for the Real Analysis courses, particularly at the University of PGRI, Palembang, Indonesia. Apart from implementing excellent teaching methods, we strongly believe that the curriculum for these courses should be adapted to the characteristics of students who are specialized in Mathematics Education and their careers after completing their degree. Although the materials in the Real Analysis will never be given to secondary school students, the participants in these classes are trained to think critically. In turn, the ability of this critical thinking is very beneficial for a good teacher. For future research, it is important to investigate the significant value of the students' responses. This investigation should involve the quantitative calculation and validation test.

\end{document}